\documentclass[11pt]{amsart}
\usepackage{amssymb}
\usepackage[english]{babel} 
\def\R{\mathbb R}

\newtheorem{theoreme}{Theorem}

\newtheorem{lemma}{Lemma}
\newtheorem{remark}{Remark}

\numberwithin{equation}{section}
\numberwithin{proposition}{section}

\begin{document}
\title[On the Benjamin-Ono equation]
{Nonlinear wave interactions for the Benjamin-Ono equation  }
\author{H. Koch}
\address{Universit\"at Dortmund, 44221 Dortmund}
\email{herbert.koch@mathematik.uni-dortmund.de}
\author{N. Tzvetkov}
\address{
D\'epartement de Math\'ematiques, Universit\'e Lille I, 59 655 Villeneuve d'Ascq Cedex
}
\email{nikolay.tzvetkov@math.univ-lille1.fr}
\selectlanguage{english}
\begin{abstract}
{
We study the interaction of suitable small and high frequency
waves evolving by the flow of the Benjamin-Ono equation.
As a consequence, we prove that the flow map of the Benjamin-Ono equation 
can not be uniformly continuous on
bounded sets of $H^s(\R)$ for $s>0$.  
} 
\end{abstract} 
\maketitle 
\section{Introduction}
In this paper we continue our investigations around the Cauchy problem for the Benjamin-Ono equation
\begin{equation}\label{BO} u_{t}+Hu_{xx}+uu_{x}=0, \end{equation}
where $H$ denotes the Hilbert transform. 
The Benjamin-Ono equation is a model for the propagation of one dimensional 
internal waves (see \cite{B}). In \cite{Tao} it was shown  
that the Benjamin-Ono equation is globaly well-posed in $H^1(\mathbb{R})$,  
(see also \cite{KK,KT1,Ponce,Saut}). 
Most likely the lower bound on $s$ is  not optimal.
The local well-posedness implies that the flow map
is continuous on $H^s(\R)$. The main purpose of this paper is to show that no
further regularity holds. 
More precisely we are going to show that the flow map of the  
Benjamin-Ono equation can not be uniformly continuous
on bounded sets of $H^s(\R)$ for $s>0$ and thus we significantly extend 
the recent work \cite{MST} of Molinet, Saut and the second 
author where it is shown   that the flow map 
of the Benjamin-Ono equation can not be of class $C^2$ as map 
on $H^s(\R)$. 

Our method of proof is based on the description of the effect of a suitable 
small low frequency perturbation to a high frequency wave evolving by the flow of the Benjamin-Ono equation.  
We believe that these considerations are of  independent interest and we refer to the next sections 
for more details. We now state our result concerning the lack of uniform continuity of the  
flow map of the Benjamin-Ono equation.
\begin{theoreme}\label{th1}
Let $s>0$. There exist two positive constants $c$ and $C$ and 
two sequences $(u_n)$ and $(\widetilde{u}_{n})$ of solutions of (\ref{BO})  
such that for every $t\in [0,1]$,
$$
\sup_{n}\|u_{n}(t,\cdot)\|_{H^{s}(\R)}+\sup_{n}\|\widetilde{u}_{n}(t,\cdot)\|_{H^{s}(\R)}\leq C\, ,
$$
$(u_n)$ and $(\widetilde{u}_{n})$ satisfy initially
$$
\lim_{n\rightarrow\infty}\|u_{n}(0,\cdot)-\widetilde{u}_{n}(0,\cdot)\|_{H^{s}(\R)}=0,
$$
but, for every $t\in [0,1]$,
\begin{equation*}
\liminf_{n\rightarrow \infty}\|u_{n}(t,\cdot)-\widetilde{u}_{n}(t,\cdot)\|_{H^{s}(\R)}\geq c\,\sin t\, .
\end{equation*}
\end{theoreme}
The case $s<-1/2$ can be easily treated
by using the high speed limits of solitary waves (see \cite{BiLi}).
The construction of the solutions of Theorem \ref{th1} relies on a separation 
of the transport and the dispersion effects. This phenomenon
becomes much less clear when $s \le 0$. The statement of failure of 
uniform continuity 
is fairly strong: None of the maps $u_0 \to u(t)$ can be uniformly continuous 
in balls. We can obtain a similar statement for lower values of $s$ 
if consider uniform continuity of the map
\begin{equation}\label{map}
H^s \ni u_0 \to u \in C([0,1];H^s)
\end{equation}
instead. For instance, in \cite{CCT2} failure of the uniform continuity 
of the map (\ref{map})
in the context of supercritical wave a Schr\"odinger equations is 
obtained but
no result for the map  $u_0 \to u(t)$,  $t>0$ is available.

We are inspired by the following observation on the
Burgers equation. If $v$ solves the Burgers equation 
$v_t = v v_x $ then so does 
$$
w(t,x) :=  v(t,x+\omega t) + \omega.
$$
However, the shift in the spatial variable is not norm continuous, considered as a map from $\R$ to the bounded
operators on $H^s(\R)$. 
Our strategy will be to construct solutions of 'fixed' frequency, localized 
in $x$, to which we add a smooth solution, which changes the speed of the wave. 

It is worth noticing that the ``instability property'' 
of the flow of the BO equation displayed by Theorem \ref{th1}
is not shared by the KdV equation which is another important model 
for the propagation of one directional waves. 
Due to the higher speed of propagation  of the linearized KdV waves, 
the flow of the KdV equation is Lipschitz continuous on bounded sets of $H^s$, $s>-3/4$
(see \cite{Bo1,Bo2,KPV1}) and therefore the local flow of the Benjamin-Ono equation
turns out to be quite different from that of the KdV equation.
Since a Picard iteration scheme implies smooth dependence on the initial data, 
it is a consequence of Theorem \ref{th1}  (or of the results of \cite{MST})
that the Picard iteration scheme can not be used to construct 
solutions to the Benjamin-Ono equation for initial data in $H^s(\R)$.

Our method to prove Theorem \ref{th1} is quite general and can be applied to many 
others PDE's (see Remark \ref{r1} below).  In particular we use only the first few integrals of the 
Benjamin-Ono and we do this  only for the range $0<s\le 3/2$. 

There is a number of recent papers dealing with the lack of uniform continuity
for nonlinear PDE's (see e.g. \cite{BGT,BGT3,CCT,KPV2}). These results are however of different nature
comparing to Theorem \ref{th1}. First, to our knowledge, in none of the cases
considered in \cite{BGT,BGT3,CCT,KPV2} local well-posedness is known. Moreover, essentially 
all papers \cite{BGT,BGT3,CCT,KPV2} deal with the self interaction of a single high frequency wave which
naturally lead to upper bounds on the Sobolev regularity $s$ which is 
consistent with the smooth dependence on the 
initial data for large $s$ holding
in the cases 
considered in \cite{BGT,BGT3,CCT,KPV2}.

The rest of the paper is organized as follows. The next section is devoted to some preliminaries.
In Sections 3 and 4, we construct families of approximate solutions of the Benjamin-Ono equation.
In Sections 5 and 6 we construct solutions which are 
close to the approximate solution.
The lack of uniform continuity of the flow map (Theorem 1), proved in Section 7 is an immediate consequence. 
\\

{\bf Notation.}  
In the sequel $\lambda$ will be a large real number, and $0<\delta<1$ 
will be fixed most of the time. For functions $f: \R \to \R$ we denote 
$ f_\lambda(x) := f(x/\lambda^{1+\delta}). $

\section{Preliminaries}
In this section, we collect several preliminary results needed for the sequel.
We first state a basic existence result for the Benjamin-Ono equation
(see \cite{ABFS,KK,Tao,KT1,Ponce,Saut}).
\begin{lemma}\label{wp}
Fix $s>3/2$ and $\sigma\in]3/2,s]$. Then for every $u_{0}\in H^{s}(\R)$ there exists a unique global
solution $u\in C(\R;H^{s}(\R))$ of (\ref{BO}) subject to the initial data $u_0$. Moreover 
$$
\|u(t,\cdot)\|_{H^{s}(\R)}\lesssim\|u_0\|_{H^{s}(\R)}
$$
if $|t|\leq \min\big(1,\, c\|u_0\|_{H^{\sigma}(\R)}^{-4}\big)$.
\end{lemma}
Even more is true: For initial data in $H^1(\R)$ there exists a unique 
global solution $u$ with $\partial_x u \in L^1_{loc}(\R; L^\infty(\R))$ see \cite{Tao}.
Nevertheless we state the  weaker existence result of Lemma \ref{wp} 
 which does not depend as heavily on the specific structure of the problem. 

The next lemma deals with the commutators of the Hilbert transform $H$  and   
scaled smooth  functions
with compact support.
\begin{lemma}\label{commutator}
Fix $0<\delta<1$ and  $\phi\in C_{0}^{\infty}(\R)$.
Then for any $N>0$ there exists a positive constant $C_N$ such that
for every $\alpha\in\R$ 
\begin{equation}\label{isa}
\left\|\left[H,\phi_\lambda \right]
\cos(\lambda x+\alpha)\right\|_{L^{2}_{x}}
\leq
C_N \,\lambda^{-N}.
\end{equation}
\end{lemma}
{\bf Proof of Lemma \ref{commutator}.}
Inequality (\ref{isa}) is equivalent to 
\begin{equation}\label{isa-bis}
\begin{split} \!\!\left\|
\int_{-\infty}^{\infty}
\frac{\phi_\lambda(x)-\phi_\lambda(y)}{x-y}
\cos(\lambda y+\alpha)dy
\right\|_{L^{2}_{x}} \hspace{-3cm}  & \\
= & \lambda^{\frac{1+\delta}2} 
\left\|  \int_{-\infty}^\infty \frac{\phi(x)-\phi(y)}{x-y} \cos(\lambda^{2+\delta} y+
\alpha) dy \right\|_{L^2} \!\! 
\\ \leq &  
C_N\, \lambda^{-N}.
\end{split} 
\end{equation}
We fix $\eta \in C^\infty_0(\R)$, which is identically $1$ on $[-1,1]$ and 
between $0$ and $1$ on $\R$.   
By writing
\[
\begin{split} 
\frac{\phi(x)-\phi(y)}{x-y}
= &
\eta(x-y) \int_{0}^{1}
\phi'(t x+(1-t)y)\, dt + (1-\eta(x-y)) \frac{\phi(x)-\phi(y)}{x-y}
\\ =: &  \beta(x,y) + \gamma(x,y) 
\end{split} 
\]
we will deduce (\ref{isa-bis}) after several integration by parts in
the $y$ variable.
  Both functions $\beta$ and $\gamma$ are smooth. Integration by parts gives
\[\begin{split} 
 \int_{-\infty}^\infty \frac{\phi(x)-\phi(y)}{x-y} \cos(\lambda^{2+\delta} y+
\alpha) dy  & \\ & \hspace{-4cm} =  
 \lambda^{-8N-4\delta} 
\int_{-\infty}^\infty( \partial_y^{4N} (\beta(x,y)+ \gamma(x,y)) ) \cos(\lambda^{2+\delta} y+ \alpha) 
dy.
\end{split} 
\]
The function $\beta$ is compactly supported and the corresponding desired
bound is obvious. 
If when differentiating $\gamma$ 
 any derivative falls on $\eta$ then  we are again in the compactly supported 
situation.  By Minkowskis inequality 
it remains to provide the uniform bounds 
\[
\begin{split}   
\int_{-\infty}^\infty \left\Vert  (1- \eta(x-y)) \partial_y^{4N}\left[ (\phi(x)-\phi(y)
)|x-y|^{-1}\right] \right\Vert_{L^2_x} dy 
& 
\\ & \hspace{-6cm}\le   \int_{-\infty}^\infty (4N)! \Vert  (1- \eta(x-y)) \phi(x)|x-y|^{-4N-1}  \Vert_{L^2_x} dy 
\\ & \hspace{-5cm} + \int_{-\infty}^\infty  \left\Vert 
 (1- \eta(x-y))\partial_y^{4N} \left[ \phi(y)|x-y|^{-1} \right]  \right\Vert_{L^2_x} dy 
\\ & \hspace{-6cm} 
\lesssim \int_{-\infty}^\infty (1- \eta(z))|z|^{-4N-1} dz  
 \Vert \phi \Vert_{L^2} 
\\ & \hspace{-5cm} + \Vert \phi \Vert_{W^{4N,1}} 
 \sup_{y} \Vert (1-\eta(x-y)) |x-y|^{-1} \Vert_{L^2}.  
\end{split} 
\]
where $\Vert \phi \Vert_{W^{4N,1}} = \sum_{j=0}^{4N} \Vert \phi^{(j)}
\Vert_{L^1} $. 
\qed

\medskip

We will complete this section by evaluating the $H^s$ norm of some 
high frequency localized smooth functions.
\begin{lemma}\label{loc}
Fix $s\geq 0$, $0<\delta<1$, $\alpha\in\R$ and  $\phi\in C_{0}^{\infty}(\R)$.
Then 
\[
\lim_{\lambda \to \infty} \lambda^{-\frac{1+\delta}{2}-s} 
\left\|\phi_\lambda(x)  \cos(\lambda x+\alpha)\right\|_{H^{s}_{x}}
= \frac1{\sqrt{2}} \Vert \phi \Vert_{L^2}. 
\]
\end{lemma}
{\bf Proof of Lemma \ref{loc}.}
Write via the rescaling $x\mapsto \lambda^{1+\delta}x$,
\begin{equation} \label{edno}\begin{split}
\Vert [(1+|D|^2)^{s/2}, \phi_\lambda]  \cos(\lambda x + \alpha)
\Vert_{L^2_x} 
= &
\\
& \hspace{-1cm}   \lambda^{\frac{1+\delta}2} 
\Vert [ (1 + \lambda^{-2-2\delta} |D|^2)^{s/2}, \phi ]
\cos(\lambda^{2+\delta} x + \alpha ) \Vert_{L^2_x}.
\end{split}
\end{equation}
Plancherel's theorem then yields   
\[ 
\begin{split} 
 \Vert [ (1 + \lambda^{-2-2\delta} |D|^2)^{s/2}, \phi ]
\exp(i\lambda^{2+\delta} x)
\Vert_{L^2_{x}}  
  \hspace{-2cm} & 
\\  = & 
\left\| 
\big[(   1+ (\lambda^{-1-\delta} \xi)^2)^{s/2}
-(   1+   \lambda^2)^{s/2}) \big] \hat \phi(\xi-\lambda^{2+\delta}) \right\|_{L^2_{\xi}}
\\ 
 = &  
\left\| \big[(   1+ (\lambda^{-1-\delta} \xi+\lambda)^2)^{s/2}
-(   1+   \lambda^2)^{s/2}) \big] \hat \phi(\xi) \right\|_{L^2_{\xi}}. 
\end{split}
\]
Since $\phi$ is a Schwartz function, the above identity and (\ref{edno}) imply  
\[ \Vert [(1+|D|^2)^{s/2}, \phi_\lambda]  \cos(\lambda x + \alpha)
\Vert_{L^2_x}  \lesssim \lambda^{\frac{1+\delta}2 + s-2-\delta}. 
\]
We have
$$
(1+ |D|^2)^{s/2} \cos(\lambda x + \alpha) = (1+ |\lambda|^2)^{s/2} 
\cos(\lambda x + \alpha).
$$
Finally, we obtain by an integration by parts
\begin{eqnarray*} 
\int_{-\infty}^\infty \Big(\phi_\lambda(x)\cos(\lambda x+\alpha)\Big)^2 \, dx 
& = & 
\int_{-\infty}^\infty \phi^2_\lambda(x)  ( \frac12 + \frac12\cos( 2\lambda x+ 2\alpha)) \, dx 
\\ & = & \frac12 \lambda^{1+\delta} \Vert \phi \Vert_{L^2}^2 + 
{\mathcal O}(\lambda^{-1}\,
\Vert \phi \Vert_{L^2_x} 
\Vert \phi^\prime \Vert_{L^2_x} 
) 
\end{eqnarray*}  
which completes the proof of Lemma \ref{loc}.  \qed

\section{First construction of approximate solutions}
In the rest of the paper, we shall use make of two smooth 
characteristic functions 
$\phi$ and $\tilde{\phi}$ in the usual manner.
Namely, let $\phi\in C_{0}^{\infty}(\R)$ be such that 
$$
\phi(x)=
\left\{
\begin{array}{ll}
0,\quad |x|>2,\\
1,\quad |x|<1
\end{array}
\right.
$$
and let $\tilde{\phi}\in C_{0}^{\infty}(\R)$ be equal to one on the 
support of $\phi$. 

For $\lambda\geq 1$, $0<\delta<1$ and $\omega\in\R$, we set
$$
U_{\lambda,\omega}(x):= 
-\omega\,\lambda^{-1}\tilde{\phi}_\lambda(x)
$$
which corresponds to the low frequency part of the approximate solution.

Next, we define the phase
$$
\Phi:=-\lambda^2t+\lambda x+  \omega t 
$$
which describes the   phase shift compared  to  linear
Benjamin-Ono waves.

Further, we define the high frequency part of the approximate solution
$$
u_{h}(t,x):= -\lambda^{-\frac{1+\delta}{2}-s}\phi_{\lambda}(x)
\cos \Phi.
$$
The next lemma states that $u_{ap}(t,x)$ defined by
\begin{equation}
u_{ap}(t,x):= U_{\lambda,\omega}(x) 
+ u_{h}(t,x)
\label{ap1}
\end{equation}
almost solves the Benjamin-Ono 
equation when $s$ is not too large and $\lambda\gg 1$. 
\begin{lemma}\label{l1}
Set
\begin{equation}
F:=(\partial_{t}+H\partial_{x}^{2})u_{ap}+
u_{ap}\,\partial_{x}u_{ap}.
\end{equation}
Then there exists a positive constant $C$ such that
for $t\in\R$, $0<\delta<1$ one has
$$
\|F(t,\cdot)\|_{L^{2}(\R)}
\leq 
C(\lambda^{-\delta-s}+\lambda^{\frac{1-\delta}{2}-2s}+ 
\omega^2 \lambda^{- \frac{5+\delta}2} 
+
|\omega|\, \lambda^{- \frac{5+3\delta}{2}} 
+
|\omega|\, \lambda^{- 2-\delta-s} 
).
$$
\end{lemma}
{\bf Proof of Lemma \ref{l1}.}
Let $L: = \partial_t + H \partial_{x}^2$. 
We compute 
\begin{eqnarray*} 
Lu_{ap}+u_{ap}\,\partial_x u_{ap}
& = &   
L U_{\lambda,\omega} +  U_{\lambda,\omega}\partial_x
U_{\lambda,\omega}
+\partial_x (  U_{\lambda,\omega} u_{h} ) 
+  u_{h}\partial_x u_{h}  + L u_{h} 
\\  
& = &   
L U_{\lambda,\omega} +  U_{\lambda,\omega}\partial_x
U_{\lambda,\omega}
\\ & & - \cos\Phi\,\, \partial_x \Big\{ U_{\lambda,\omega} 
\lambda^{-\frac{1+\delta}{2}-s}\phi_\lambda \Big\}
\\
& &
+  u_{h} \partial_x u_{h} 
\\ & & 
 - \lambda^{-\frac{1+\delta}{2}-s} 
\Big[H \partial_x^2 ,\phi_\lambda  \Big] \cos \Phi
\\ & & - \lambda^{-\frac{1+\delta}{2}-s}\phi_\lambda 
(L+U_{\lambda,\omega}\partial_x)
\cos \Phi 
\\   
& =: & F_1 + F_2 + F_3 + F_4+F_5,  
\end{eqnarray*} 
where we put the part of 
$\partial_x( U_{\lambda,\omega} u_{h})$ with the
derivative on $\cos\Phi$ into $F_5$.   
Using that $\phi_{\lambda}\tilde{\phi}_\lambda=\phi_{\lambda}$,
we readily obtain that $F_5$ vanishes which is the crucial cancelation.
It is worth noticing that both
\begin{eqnarray}\label{I}
- \lambda^{-\frac{1+\delta}{2}-s}\phi_\lambda\, L(\cos \Phi) 
\end{eqnarray}
and
\begin{eqnarray}\label{II}
-\lambda^{-\frac{1+\delta}{2}-s}\phi_\lambda\, U_{\lambda,\omega}\,(\partial_x \cos\Phi)
\end{eqnarray}
are ``big'' in $L^2$ compared to $u_{ap}$.
Notice also that the term (\ref{I}) comes from the linear part while (\ref{II}) is a contribution coming
from the nonlinear term.

Next we expand $F_4$ as 
\[
F_4 = 
\lambda^{\frac{3-\delta}2-s}  [H, \phi_{\lambda}] \cos\Phi
\,+\,2\lambda^{ - \frac{1+\delta}2 - s -\delta}  
H \left\{(\phi^{\prime})_\lambda   \sin  \Phi\right\} 
- \lambda^{-\frac52 -\frac{5\delta}2 - s } 
H \left\{ (\phi^{\prime\prime})_\lambda  \cos\Phi \right\}. 
\]
The first term is controlled in $L^2$ by Lemma \ref{commutator}. 
The $L^2$ norm of the other terms 
are readily estimated by $ c \lambda^{ - \delta - s}$. 
The $L^2$ norm of $F_3$ is easily controlled by 
$ c \lambda^{ \frac{1-\delta}2 - 2s } $. 
It remains 
to estimate $F_1$, $F_2$ and $F_3$. The bound for the latter is obvious. 
The term $F_1$ can be handled as follows
$$
\|F_1\|_{L^2}
\lesssim
\Vert U_{\lambda,\omega}  \partial_x U_{\lambda,\omega}
\Vert_{L^2} +  
\Vert H U_{\lambda,\omega}^{\prime\prime}  \Vert_{L^2} 
\lesssim  
\omega^2 \lambda^{-\frac{5+\delta}2}
+
|\omega| \lambda^{-\frac{5+3\delta}2}.
$$
The term $F_2$ is estimated as
$$ 
\Vert F_2(t,\cdot) \Vert_{L^2}  
\lesssim
|\omega|\lambda^{-2-\delta-s} 
$$ 
This completes the proof of Lemma \ref{l1}.
\qed
\begin{remark}
Notice that if $0<s<2$ the bound of Lemma \ref{l1} on $F$ can be 
simply written as
\begin{eqnarray}\label{commnet}
\|F(t,\cdot)\|_{L^{2}(\R)}
\lesssim \lambda^{-\frac{\min\{\delta,1-\delta\}}{2}-s}
\end{eqnarray}
if $|\omega| \leq 1$.
The bound (\ref{commnet}) implies  that at least in
$L^{2}(\R)$, $u_{ap}(t,x)$ is a good approximate solution of the
Benjamin-Ono equation.
\end{remark} 

\section{Refined approximate solutions}
When $s\geq 2$ one has to  modify slightly the construction 
of approximate solutions, presented in the previous section.
To avoid small frequency residual terms, 
we will chose the small frequency part of $u_{ap}$ to be a solution 
of the Benjamin-Ono equation with small frequency initial data.
Let $u_{low}(t,x)$ be the solution of (\ref{BO}) with initial data
\begin{equation}\label{low}
u_{low}(0,x)=-\omega\lambda^{-1}\tilde{\phi}_\lambda(x),\quad 0<\delta<1,
\quad \omega\in\R.
\end{equation}
In the next lemma, we collect several bounds for $u_{low}(t,x)$.
\begin{lemma}\label{ulow}  
Let $k\geq 0$. Then the following estimates hold, if $|t|\leq 1$, $\lambda\gg 1$ and 
$|\omega|\ll\lambda^{\frac{1-\delta}{2}}$,
\begin{eqnarray}
\|\partial_x^k u_{low}(t,\cdot)\|_{L^2(\R)} 
& \lesssim & |\omega|\,\lambda^{-\frac{1-\delta}{2}-k(1+\delta)},
\label{parvo}
\end{eqnarray}
\begin{eqnarray}
\|\partial_x u_{low}(t,\cdot) \|_{L^{\infty}(\R)} & \lesssim & |\omega|\, \lambda^{-2-\delta},
\label{vtoro}
\end{eqnarray}
\begin{eqnarray}
\| u_{low}(t,\cdot)- u_{low}(0,\cdot) \|_{L^2(\R)} & \lesssim & |\omega|\, \lambda^{-2-\delta}.
\label{treto}
\end{eqnarray}
\end{lemma} 
{\bf Proof of Lemma \ref{ulow}.}
Rescale by setting 
\begin{eqnarray}\label{v}
v(t,x) :=  \lambda^{1+\delta} u_{low}(\lambda^{2+2\delta } t, \lambda^{1+\delta} x).
\end{eqnarray}
Then $v$ is again a solution of the Benjamin-Ono equation.
Since $v(0,x)= -\omega \lambda^{\delta}\tilde \phi(x)$ we readily obtain  for any $s\geq 0$ the bound
$$
\|v(0,\cdot)\|_{H^s}
\lesssim
|\omega|\,\lambda^{\delta}
$$
and therefore by Lemma \ref{wp}
\begin{eqnarray}\label{scaling-bis}
\|v(t,\cdot)\|_{H^s}\lesssim |\omega|\,\lambda^{\delta},
\end{eqnarray}
if $|t|\leq \min(1,|\omega|^{-4}\lambda^{-4\delta})$
and $s>3/2$. 
But since the right hand-side of (\ref{scaling-bis}) contains a constant 
which is uniformly bounded for bounded  $s$, we conclude that (\ref{scaling-bis}) is valid for any real $s$.
The Sobolev embedding and (\ref{scaling-bis}) now give
\begin{eqnarray}\label{scaling}
\|v_{x}(t,\cdot)\|_{L^{\infty}}\lesssim |\omega|\,\lambda^{\delta},
\end{eqnarray}
if $|t|\leq\min(1,|\omega|^{-4}\lambda^{-4\delta})$.

Using (\ref{v}), we deduce from (\ref{scaling}) by scaling back that
\begin{eqnarray}\label{lll}
\|\partial_{x}u_{low}(t,\cdot)\|_{L^{\infty}}
\lesssim
|\omega|\,\lambda^{-2-\delta}\, ,
\end{eqnarray}
if $|t|\leq 1$ which proves (\ref{vtoro}).

We now turn to the proof of (\ref{parvo}) and (\ref{treto}).
Differentiating (\ref{v}) and using (\ref{scaling-bis}) (with $s=k$) yields
\begin{eqnarray}\label{llll}
\|\partial_{x}^{k}u_{low}(t,\cdot)\|_{L^2}
\lesssim
|\omega|\,\lambda^{-\frac{1-\delta}{2}-k(1+\delta)},\quad k=0,1,2,\dots
\end{eqnarray}
if $|t|\leq 1$.
Estimate (\ref{llll}) proves (\ref{parvo}). Next,
using (\ref{lll}), (\ref{llll}) and the equation satisfied  by $u_{low}$ gives 
$$
\|\partial_{t}u_{low}(t,\cdot)\|_{L^2}
\lesssim
\|\partial_{x}^{2}u_{low}(t,\cdot)\|_{L^2}
+
\|\partial_{x}u_{low}(t,\cdot)\|_{L^{\infty}}
\|u_{low}(t,\cdot)\|_{L^2}
\lesssim
|\omega|\,\lambda^{-2-\delta},
$$
if $|t|\leq 1$.
We now observe that (\ref{treto}) can be deduced from the above bound
via the fundamental theorem of calculus, applied to $u_{low}$ in the time variable.
This completes the proof of Lemma \ref{ulow}. 
\qed 
\\ 

We now set for $\lambda\geq 1$, $0<\delta<1$ and $|\omega| \ll 
\lambda^{\frac{1-\delta}2}$,
\begin{equation}
u_{ap}(t,x):=
u_{low}(t,x)-
\lambda^{-\frac{1}{2}-\frac{\delta}{2}-s}\phi_\lambda(x)
\cos(-\lambda^2t+\lambda x-\lambda\, t\, u_{low}(0,x)).
\label{ap2}
\end{equation}
The above function is an approximate solution for
$\lambda\gg 1$ and $s>0$.
\begin{lemma}\label{l2}
Let $s>0$, $0<\delta<1$,  $|\omega| \ll
\lambda^{\frac{1-\delta}2}$ and $|t|\le 1$. Set
$$
F:=(\partial_{t}+H\partial_{x}^{2})u_{ap}+u_{ap}\,\partial_{x}u_{ap}.
$$
Then there exist positive constants 
$C$ and $\lambda_0$ such that for 
$\lambda\geq \lambda_0$
one has
$$
\|F(t,\cdot)\|_{L^{2}(\R)}
\leq
C\left( 
\lambda^{-\delta-s} + \lambda^{\frac{1-\delta}2-2s}
\right).
$$ 
\end{lemma}
{\bf Proof of Lemma \ref{l2}.}
Set $\Phi:=-\lambda^2t+\lambda x+\omega\, t $.
We observe that
\begin{eqnarray*}
u_{ap}(t,x)=u_{low}(t,x)-
\lambda^{-\frac{1}{2}-\frac{\delta}{2}-s}\phi_{\lambda}(x)\cos \Phi.
\end{eqnarray*}
Then, as in the proof of Lemma \ref{l1}, 
we can write again  
$$
(\partial_{t}+H\partial_{x}^{2})u_{ap}+u_{ap}\,\partial_{x}u_{ap}
= F_1 + F_2 + F_3 + F_4 +F_5,
$$
where $U_{\lambda,\omega}$ is simply replaced by $u_{low}$.
Since $u_{low}$ is a solution of the Benjamin-Ono equation, we deduce
that $F_1=0$ which eliminates the problem with small frequency residual terms.
A difficulty however appears since now $F_5$ is not vanishing anymore.
We will however be able to control $F_5$ by the aid of Lemma \ref{ulow}.  
Using that
$\phi_{\lambda}\tilde{\phi}_\lambda=\phi_{\lambda}$, we readily obtain
that
$$
F_5  = 
\lambda^{\frac{1-\delta}{2}-s} (u_{low}(t,x) - u_{low}(0,x) )
\phi_\lambda(x)
\sin \Phi .
$$
Using Lemma \ref{ulow}, we get
\begin{eqnarray}\label{F_5}
\|F_5(t,\cdot)\|_{L^2} \lesssim 
\lambda^{\frac{1-\delta}{2}-s}\,
|\omega| \lambda^{-2-\delta}
\lesssim
\lambda^{-1-2\delta-s}.
\end{eqnarray}
It remains to bound $F_2$, $F_3$ and $F_4$.
Observe that $F_3$ and $F_4$ are exactly as in Lemma \ref{l1}
and therefore
\begin{eqnarray}\label{F34}
\|F_3(t,\cdot)\|_{L^2} +\|F_4(t,\cdot)\|_{L^2} 
\lesssim \lambda^{-\delta-s}+\lambda^{\frac{1-\delta}{2}-2s}.
\end{eqnarray}
The term $F_2$ reads
$$
F_2  = -\cos\Phi\,\,\partial_x \Big\{u_{low}(t,x) 
\lambda^{-\frac{1+\delta}{2}-s}\phi_\lambda(x)\Big\}.
$$
Using Lemma \ref{ulow} and the assumption on $|\omega|$, we obtain
\begin{multline}\label{F2}
\|F_2(t,\cdot)\|_{L^2}\lesssim 
\lambda^{-s}\|\partial_x u_{low}(t,\cdot)\|_{L^{\infty}}
+
\lambda^{-\frac{3+3\delta}{2}-s}
\|u_{low}(t,\cdot)\|_{L^{2}}
\lesssim
\lambda^{-\frac{3+3\delta}{2}-s}.
\end{multline}
Collecting (\ref{F_5}), (\ref{F34}) and (\ref{F2}) 
completes the proof of Lemma \ref{l2}.
\qed
\section{Analysis of the difference equation}
In order to prove Theorem \ref{th1}, we need to show that
the family of approximate solutions constructed in Sections 3 and 4 are
indeed close to the ``real'' solutions of the Benjamin-Ono equation
at least up to $t=1$.
\begin{theoreme}\label{th2}
Let  $1-s<\delta< 1$ and $|\omega| \ll \lambda^{\frac{1-\delta}2}$. Let $u_{\omega,\lambda}$ 
be the unique global solution of the Benjamin-Ono equation subject to initial data
$$
u_{\omega,\lambda}(0,x)
=
-\omega\,\lambda^{-1}\tilde{\phi}_\lambda(x)
-
\lambda^{-\frac{1}{2}-\frac{\delta}{2}-s} \phi_\lambda(x)
\cos \lambda x.
$$
Then the identity
$$
u_{\omega,\lambda}(t,x) = -\lambda^{-\frac{1}{2}-\frac{\delta}{2}-s} 
\phi_\lambda(x)
\cos (-\lambda^2 t + \lambda x + \omega t)  + 
{\mathcal O}\left(\lambda^{- \frac{\min\{\delta, 1-\delta\}}{4(s+2)} } 
+ |\omega|\, \lambda^{-\frac{1-\delta}2}\right)
$$
holds in $H^{s}_{x}(\R)$, uniformly in $t\in[0,1]$ and $\lambda\gg 1$.
\end{theoreme}
\begin{remark}
The theorem  describes  the short time nonlinear interaction
between some low and high frequency waves.
If $\omega=0$, the approximate solution propagates as a high frequency
linear Benjamin-Ono wave.
When $\omega\neq  0$, the approximate solution propagates 
as a high frequency linear dispersive wave with modified propagation speed
which is the crucial nonlinear effect.
\end{remark}
{\bf Proof of Theorem \ref{th2}.}
The first step is to bound $u_{\omega,\lambda}$ in high Sobolev norms. \\
Let $s>\frac{3}{2}$. 
Observe that for $\frac{3}{2}<\sigma<s$
$$
\|u_{\omega,\lambda}(0,\cdot)\|_{H^{\sigma}}
\lesssim
\lambda^{\sigma-s}+|\omega|\lambda^{-\frac{1-\delta}{2}}.
$$
Therefore for $k\geq s$, it follows from Lemma \ref{wp} that
\begin{equation}\label{high1} 
\Vert u_{\omega,\lambda}(t,\cdot) \Vert_{H^k} 
\lesssim \Vert u_{\omega,\lambda}(0,\cdot) \Vert_{H^k}
\lesssim \lambda^{k-s}, \quad |t|\leq 1,\quad \lambda\gg 1.
\end{equation}
Let $0<s\le \frac32$.
Using the conservation laws associated to the Benjamin-Ono equation
(see Lemma 3.3.2 of \cite{ABFS}), we  get the following bound uniformly 
in $t\in\R$ 
\begin{eqnarray}\label{CL}
\|u_{\omega,\lambda}(t,\cdot)\|_{H^{2}}
\lesssim
\|u_{\omega,\lambda}(0,\cdot)\|_{H^{2}}
+
\|u_{\omega,\lambda}(0,\cdot)\|_{L^{2}}^{5}
\lesssim
1+  \lambda^{2-s},
\end{eqnarray}
and therefore we obtain 
\begin{equation} \label{high2}
\|u_{\omega,\lambda}(t,\cdot)\|_{H^{2}}
\lesssim
\lambda^{2-s},\quad t\in\R.
\end{equation}
Let $u_{ap}$ as in \eqref{ap2} and
$
v_{\omega,\lambda}:=u_{\omega,\lambda}-u_{ap}.
$
The aim is to show that $v_{\omega,\lambda}$ is small comparing to $u_{ap}$ in the $H^s$ norm.

Due to Lemma \ref{ulow}, we get
$$
\Vert u_{low}(t,\cdot)\Vert_{H^s} 
\lesssim |\omega|\, \lambda^{-\frac{1-\delta}2},
$$
if $|t|\leq 1$.
Next, using Lemma \ref{loc}, we obtain the bound
$$
\|u_{ap}(t,\cdot)\|_{H^{k}(\R)}
\lesssim \lambda^{k-s},
$$
if $|t|\leq 1$ and $k\geq s$.

Therefore using (\ref{high1}) and (\ref{high2}), we get the bounds for the  
high Sobolev norms 
\begin{eqnarray} 
\label{firstbound}
\|v_{\omega,\lambda}(t,\cdot)\|_{H^{k}}
\lesssim
\lambda^{k-s}, 
\end{eqnarray}
if $|t|\leq 1$ and $\frac{3}{2}<s<k$, and
\begin{eqnarray} 
\label{firstbound-bis}
\|v_{\omega,\lambda}(t,\cdot)\|_{H^{2}}
\lesssim
\lambda^{2-s}, 
\end{eqnarray}
if $t\in\R$ and $0<s<\frac{3}{2}$.

The second step provides a good bound of the $L^2$ norm of 
$v_{\omega,\lambda}$. Clearly
\begin{equation}\label{Eq-v}
\left\{
\begin{array}{l}
(\partial_{t}+H\partial^{2}_{x})v_{\omega,\lambda}
+
v_{\omega,\lambda}\,\partial_{x}v_{\omega,\lambda}
+
\partial_{x}(u_{ap}\,v_{\omega,\lambda})+
F=0,
\\
v_{\omega,\lambda}(0,x)=0
\end{array}
\right.
\end{equation}
with 
$$ F = (\partial_t + H \partial_x^2 ) u_{ap} +u_{ap}\partial_x u_{ap}, $$
which satisfies 
$$\|F(t,\cdot)\|_{L^{2}(\R)} \lesssim \lambda^{- \frac{\min \{ \delta, 1-\delta\}}2 -s}$$
by Lemma \ref{l2} and the assumption $1-s<\delta<1$.

The second endpoint is the $L^2$ estimate 
\begin{equation}\label{second}  
\Vert v_{\omega,\lambda}(t,\cdot) \Vert_{L^2} \lesssim \lambda^{ - \frac{\min \{ \delta, 1-\delta\}}2 -s },\quad |t|\leq 1. 
\end{equation}  
To prove \eqref{second}, we multiply (\ref{Eq-v}) by $v_{\omega,\lambda}$ and integrate over the real line, 
\begin{eqnarray*}
\frac{d}{dt}\|v_{\omega,\lambda}(t,\cdot)\|_{L^2}^{2}
\lesssim
\|\partial_{x}u_{ap}(t,\cdot)\|_{L^{\infty}}
\|v_{\omega,\lambda}(t,\cdot)\|_{L^2}^{2}
+
\|v_{\omega,\lambda}(t,\cdot)\|_{L^2}
\|F (t,\cdot)\|_{L^{2}}
\end{eqnarray*}
hence, since we have for $1-s <\delta<1$ 
\[ \|\partial_x u_{ap}(t,\cdot)\|_{L^{\infty}}\lesssim 
\|\partial_x u_{low}(t,\cdot)\|_{L^{\infty}}
+
\lambda^{\frac{1-\delta}2-s} 
\lesssim
|\omega|\lambda^{-2-\delta}
+
\lambda^{\frac{1-\delta}2-s} 
\ll 1,
\] 
we readily get the bound \eqref{second}.

We now complete the proof by an interpolation argument.
Let first $s>\frac{3}{2}$. 
Choose $k\in [s+\frac{1}{2},s+2]$ and interpolate between
\eqref{firstbound} and \eqref{second} as follows
$$
\Vert v_{\omega,\lambda}(t,\cdot) \Vert_{H^s} 
\le 
\Vert v_{\omega,\lambda}(t,\cdot) \Vert_{L^2}^{\frac{k-s}k} 
\Vert v_{\omega,\lambda}(t,\cdot) \Vert_{H^k}^{\frac{s}k}
\lesssim\lambda^{ -\frac{\min\{\delta,1-\delta\}}{4(s+2)}}.   
$$
If $s\le \frac32$ we obtain the same estimate by using $k=2$ in the
interpolation and (\ref{firstbound-bis}) instead of \eqref{firstbound}.
This completes the proof of Theorem \ref{th2}. 
\qed
\begin{remark}\label{r1}
Notice that to derive estimate \eqref{CL} one needs to exploit
the higher conservation laws for the Benjamin-Ono equation.
This fact permits us to have an ansatz for the solution up to time
one as claimed in Theorem \ref{th2}.
If $s > 3/2$ we use Lemma \ref{wp} instead. 

The method of proof can be generalized to many other equations.
For example the corresponding to Theorem \ref{th2} 
result in the context of the KdV equation provides 
a family of essentially linear KdV waves 
($\omega\rightarrow 0$) with the 
same initial data as approximate solutions and thus no instability 
property of the flow is displayed.
\end{remark}
\section{Lack of uniform continuity}
In this section, we complete the proof of Theorem \ref{th1}. We apply Theorem \ref{th2} 
with $\omega=\pm 1$ and $\lambda=2^n$
to obtain two families $(u_{1,2^n})$ and $(u_{-1,2^n})$
of solutions to  the Benjamin-Ono equation. 
Notice that
$$
\|u_{1,2^n}(0,\cdot)-u_{-1,2^n}(0,\cdot)\|_{H^{s}}
\lesssim 2^{\frac{(\delta-1)n}{2}} 
$$
and moreover due to Theorem \ref{th2}, setting  $\kappa= -2^{2n} t+ 2^n x$, we arrive at
$$
\|u_{1,2^n}(t,\cdot)-u_{-1,2^n}(t,\cdot)\|_{H^{s}}=
\Vert 
2^{-(\frac{1+\delta}2+s)n} \phi_{2^n}(x) (\cos(\kappa   +t ) - 
\cos(\kappa-t )) \Vert_{H^s_x} + o(1), 
$$
if $|t|\leq 1$ and where $o(1)\rightarrow 0$ as $n\rightarrow\infty$.
Then using Lemma \ref{loc}  we get
\begin{eqnarray*}
\Vert 
2^{-(\frac{1+\delta}2+s)n} \phi_{2^n}(x) (\cos(\kappa   +t ) - 
\cos(\kappa-t )) \Vert_{H^s_x}=\sqrt{2}|\sin t|  \,\,
\Vert \phi \Vert_{L^2}+o(1).
\end{eqnarray*} 
The proof of Theorem \ref{th1} is completed.
\qed
 
\bigskip

{\bf Acknowledgments.}
This work was initiated during a visit of the second author at University
of Dortmund. 
We acknowledge the support for this work by the European Commission
through the IHP Network HARP ``Harmonic Analysis and Related Problems".
The second author would like to thank
the Mathematics department of the 
University of Dortmund for the kind hospitality.

%\backmatter

\end{document}